\documentclass[15pt, a4paper]{article}
\usepackage{}
\usepackage{amsfonts}
\usepackage{mathbbold}
\usepackage{bbm}
\usepackage{mathrsfs}
\usepackage{latexsym}
\usepackage{graphicx}
\usepackage{amssymb}

\newtheorem{theorem}{Theorem}[section]
\newtheorem{lemma}[theorem]{Lemma}
\newtheorem{corollary}[theorem]{Corollary}

\title{{\Large \bf Some results on extremal spectral radius of hypergraph\thanks{Supported by National natural science foundation of China (NSFC)
(Nos. Nos. 12171222, 12101285, 12071411), Natural science foundation of Guangdong province (No. 2021A1515010254), Foundation of Lingnan Normal
University(ZL1923).}}}
\author{Guanglong Yu$^{a}$
~ Lin Sun$^{a}$\thanks{Corresponding authors, E-mail addresses:
yglong01@163.com (G. Yu), sunlin@lingnan.edu.cn (L. Sun).} ~  Hailiang Zhang$^b$ ~ Gang Li$^{c}$ ~
\\ ~ \\
{\footnotesize $^a$Department of Mathematics, Lingnan Normal
University,  Zhanjiang, Guangdong, 524048, P.R.China}\\
{\footnotesize $^b$Department of Mathematics, Taizhou University, Linhai, Zhejiang, 317000, P.R.China}\\
{\footnotesize $^c$College of Mathematics and Systems Science, Shandong University of Science and Technology,}\\
{\footnotesize  Qingdao, Shandong, 266590, P.R.China}}
\date{}

\begin{document}
%\openup 1.0\jot
\maketitle

\begin{abstract}
For a $hypergraph$ $\mathcal{G}=(V, E)$ with a nonempty vertex set $V=V(\mathcal{G})$ and an edge set $E=E(\mathcal{G})$, its $adjacency$ $matrix$ $\mathcal {A}_{\mathcal{G}}=[(\mathcal {A}_{\mathcal{G}})_{ij}]$ is defined as
$(\mathcal {A}_{\mathcal{G}})_{ij}=\sum_{e\in E_{ij}}\frac{1}{|e| - 1}$, where $E_{ij} = \{e\in E\, |\, i, j \in e\}$.
The $spectral$ $radius$ of a hypergraph $\mathcal{G}$, denoted by $\rho(\mathcal {G})$, is the maximum modulus among all eigenvalues of $\mathcal {A}_{\mathcal{G}}$.
In this paper, we get a formula about the spectral radius which link the ordinary graph and the hypergraph, and represent some results on the spectral radius changing under some graphic structural perturbations. Among all $k$-uniform ($k\geq 3$) unicyclic hypergraphs with fixed number of vertices, the hypergraphs with the minimum, the second the minimum spectral radius are completely determined, respectively; among all $k$-uniform ($k\geq 3$) unicyclic hypergraphs with fixed number of vertices and fixed girth, the hypergraphs with the maximum spectral radius are completely determined; among all $k$-uniform ($k\geq 3$) $octopuslike$ hypergraphs with fixed number of vertices, the hypergraphs with the minimum spectral radius are completely determined. As well, for $k$-uniform ($k\geq 3$) $lollipop$ hypergraphs, we get that the spectral radius decreases with the girth increasing.


\end{abstract}
\baselineskip 18.6pt

\section{Introduction}

\ \ \ \
In the past twenty years, different hypermatrices or tensors for hypergraphs have been developed to explore spectral hypergraph theory. Many interesting spectral properties of hypergraphs have been explored \cite{Ban.CB}-\cite{COQY}, \cite{LQIE}-\cite{JZLS}.
Recently, A. Banerjee \cite{A.Ban} introduced an adjacency matrix and use its spectrum so that some spectral and structural
properties of hypergraphs are revealed. In this paper, we go on studying the spectra of hypergraphs according to the adjacency matrix introduced in \cite{A.Ban}.

Now we recall some notations and definitions related to hypergraphs.
For a set $S$, we denote by $|S|$ its cardinality. A $hypergraph$ $\mathcal{G}=(V, E)$ consists of a nonempty vertex set $V=V(\mathcal{G})$ and an edge set $E=E(\mathcal{G})$. where each edge $e\in E(\mathcal{G})$ is a subset of $V(\mathcal{G})$ containing at least two vertices. The cardinality $n=|V(\mathcal{G})|$ is called the order; $m=|E(\mathcal{G})|$ is called the edge number of hypergraph $\mathcal{G}$. Denote by $t$-set a set with size (cardinality) $t$. We say that a hypergraph $\mathcal{G}$ is $uniform$ if its every edge has the same size, and call it $k$-$uniform$ if its every edge has size $k$ (i.e. every edge is a $k$-subset). It is known that a $2$-uniform graph is always called a ordinary graph or graph for short.

For a hypergraph $\mathcal{G}$, we define $\mathcal{G}-e$ ($\mathcal{G}+e$)
to be the hypergraph obtained from $\mathcal{G}$ by deleting the edge $e\in
 E(\mathcal{G})$ (by adding a new edge $e$ if $e\notin
 E(\mathcal{G})$); for an edge subset $B\subseteq E(\mathcal{G})$, we define $\mathcal{G}-B$
to be the hypergraph obtained from $\mathcal{G}$ by deleting each edge $e\in
 B$; for a vertex subset $S\subseteq V(\mathcal{G})$, we define $\mathcal{G}-S$ to be the hypergraph obtained from $\mathcal{G}$ by deleting all the vertices in $S$ and deleting the edges incident with any vertex in $S$. For two $k$-uniform hypergraphs $\mathcal{G}_{1}=(V_{1}, E_{1})$ and $\mathcal{G}_{2}=(V_{2}, E_{2})$, we say the two graphs are $isomorphic$ if there is a bijection $f$ from $V_{1}$ to $V_{2}$, and there is a bijection $g$ from $E_{1}$ to $E_{2}$ that maps each edge $\{v_{1}$, $v_{2}$, $\ldots$, $v_{k}\}$ to $\{f(v_{1})$, $f(v_{2})$, $\ldots$, $f(v_{k})\}$.

In a hypergraph, two vertices are said to be $adjacent$ if
both of them are contained in an edge. Two edges are said to be $adjacent$
if their intersection is not empty. An edge $e$ is said to be $incident$ with a vertex $v$ if
$v\in e$. The $neighbor$ $set$ of vertex $v$ in hypergraph $\mathcal{G}$, denoted by $N_{\mathcal{G}}(v)$, is the set of vertices adjacent to $v$ in $\mathcal{G}$. The $degree$ of a vertex $v$ in $\mathcal{G}$, denoted by $deg_{\mathcal{G}}(v)$ (or $deg(v)$ for short), is the number of the edges incident with $v$. A vertex of degree $1$ is called a $pendant$ $vertex$. A $pendant$ $edge$ is an edge with at most one vertex of degree more than one and other vertices in this edge being all pendant vertices.

In a hypergraph, a $hyperpath$ of length $q$ ($q$-$hyperpath$) is defined to be an alternating sequence
of vertices and edges $v_{1}e_{1}v_{2}e_{2}\cdots v_{q}e_{q}v_{q+1}$ such that
(1) $v_{1}$, $v_{2}$, $\ldots$, $v_{q+1}$ are all distinct vertices;
(2) $e_{1}$, $e_{2}$, $\ldots$, $e_{q}$ are all distinct edges;
(3) $v_{i}$, $v_{i+1}\in e_{i}$ for $i = 1$, $2$, $\ldots$, $q$;
(4) $e_{i}\cap e_{i+1}=v_{i+1}$ for $i = 1$, $2$, $\ldots$, $q-1$; (5) $e_{i}\cap e_{j}=\emptyset$ if $|i-j|\geq 2$.
If there is no discrimination, a hyperpath is sometimes written as $e_{1}e_{2}\cdots e_{q-1}e_{q}$, $e_{1}v_{2}e_{2}\cdots v_{q}e_{q}$ or $v_{1}e_{1}v_{2}e_{2}\cdots v_{q}e_{q}$. A $hypercycle$ of length $q$ ($q$-$hypercycle$) $v_{1}e_{1}v_{2}e_{2}\cdots v_{q-1}e_{q-1}v_{q}e_{q}v_{1}$ is obtained from a hyperpath $v_{1}e_{1}v_{2}e_{2}\cdots v_{q-1}e_{q-1}v_{q}$ by adding a new edge $e_{q}$ between $v_{1}$ and $v_{q}$ where $e_{q}\cap e_{1}= \{v_{1}\}$, $e_{q}\cap e_{q-1}= \{v_{q}\}$, $e_{q}\cap e_{j}=\emptyset$ if $j\neq1, q-1$ and $|q-j|\geq 2$. The length of a hyperpath $P$ (or a hypercycle $C$), denoted by $L(P)$ (or $L(C)$), is the number of the edges in $P$ (or $C$). In a hypergraph $\mathcal{G}$, the girth of $\mathcal{G}$, denoted by $g$ or $g_{\mathcal{G}}$, is the length of the shortest cycle in $\mathcal{G}$. A hypergraph $\mathcal{G}$
is connected if there exists a hyperpath from $v$ to $u$ for all $v, u \in V$,
and $\mathcal{G}$ is called $acyclic$ if it contains no hypercycle.

Recall that a tree is an ordinary graph which is 2-uniform, connected and acyclic. A $supertree$ is similarly defined to be a hypergraph which is both connected and acyclic.
Clearly, in a supertree, its each pair of the
edges have at most one common vertex. Therefore, the edge number of a $k$-uniform supertree of order $n$ is $m=\frac{n-1}{k-1}$. A connected $k$-uniform hypergraph with $n$ vertices and $m$ edges is $r$-$cyclic$ if $n-1 =(k-1)m-r$. For $r=1$, it is called a $k$-uniform $unicyclic$ $hypergraph$; for $r=0$, it is a $k$-uniform supertree. In \cite{COQY}, C. Ouyang et al. proved that a simple connected $k$-graph $G$ is unicyclic (1-cyclic) if and only if it has only one cycle. From this, for a unicyclic hypergraph $G$ with unique cycle $C$, it follows that (1) if $L(C)=2$, then the two edges in $C$ have exactly two common vertices, and $|e\cap f|\leq 1$ for any two edges $e$ and $f$ not in $C$ simultaneously; (2) if $L(C)\geq 3$, then any two edges in $G$ have at most one common vertices; (3) every connected component of $G-E(C)$ is a supertree.

Denote by $\mathscr{C}(n, k)=v_{1}e_{1}v_{2}e_{2}\cdots v_{m-1}e_{m-1}v_{m}e_{m}v_{1}$ the $k$-uniform hypercycle of order $n$ with $m=\frac{n}{k-1}$.

Let $\mathcal{C}_{g}=v_{1}e_{1}v_{2}e_{2}\cdots v_{q-1}e_{q-1}v_{g}e_{g}v_{1}$ be the $k$-uniform $g$-hypercycle with $e_{i}=\{v_{i}$, $v_{a(i,1)}$, $v_{a(i,2)}$, $\ldots$, $v_{a(i,k-2)}$, $v_{i+1}\}$ for $i=1, 2, \ldots, g-1$, $e_{g}=\{v_{g}$, $v_{a(g,1)}$, $v_{a(g,2)}$, $\ldots$, $v_{a(g,k-2)}$, $v_{1}\}$. We denote by $\mathscr{U}^{\ast}(n,k,g)$ the $k$-uniform unicyclic hypergraph of order $n$ obtained by attaching $\frac{n}{k-1}-g$ pendant edge at vertex $v_{1}$ of $\mathcal{C}_{g}$. Let $Lop_{g,k,s;1}$ be the first type $k$-uniform $lollipop$ hypergraph, which is obtained by attaching a pendant hyperpath of length $s$ at vertex $v_{1}$ of $\mathcal{C}_{g}$; $Lop_{g,k,s;2}$ be the second type $k$-uniform $lollipop$ hypergraph, which is obtained by attaching a pendant hyperpath of length $s$ at vertex $v_{a(1,1)}$ of $\mathcal{C}_{g}$ (see Fig. 1.1). A $k$-uniform ($k\geq 3$) $octopuslike$ hypergraph of order $n$, denoted by $\mathcal{O}_{n,k,g}$, is a unicyclic hypergraph obtained by attaching a supertree at a vertex of $\mathcal{C}_{g}$. Denote by $\mathcal{O}_{n,k,g;1}$ a $k$-uniform octopuslike hypergraph of order $n$, which is a unicyclic hypergraph obtained by attaching a supertree $T_{1}$ at vertex $v_{1}$ of $\mathcal{C}_{g}$; denote by $\mathcal{O}_{n,k,g;2}$ a $k$-uniform octopuslike hypergraph of order $n$, which is a unicyclic hypergraph obtained by attaching a supertree $T_{2}$ at $v_{a(1,1)}$ of $\mathcal{C}_{g}$ (see Fig. 1.1). Obviously, $Lop_{g,k,s;1}$ of order $n$ is a special $\mathcal{O}_{n,k,g;1}$; $Lop_{g,k,s;2}$ of order $n$ is a special $\mathcal{O}_{n,k,g;2}$.

\setlength{\unitlength}{0.65pt}
\begin{center}
\begin{picture}(672,383)
\put(50,314){\circle*{4}}
\put(50,300){\circle*{4}}
\put(49,330){\circle*{4}}
\put(97,366){\circle*{4}}
\qbezier(49,330)(58,361)(97,366)
\qbezier(49,330)(87,334)(97,366)
\put(150,350){\circle*{4}}
\qbezier(97,366)(131,372)(150,350)
\qbezier(97,366)(118,338)(150,350)
\put(169,299){\circle*{4}}
\qbezier(150,350)(175,338)(169,299)
\qbezier(150,350)(138,319)(169,299)
\put(148,257){\circle*{4}}
\qbezier(169,299)(173,265)(148,257)
\qbezier(169,299)(142,283)(148,257)
\put(98,244){\circle*{4}}
\qbezier(148,257)(127,239)(98,244)
\qbezier(98,244)(116,268)(148,257)
\put(50,269){\circle*{4}}
\qbezier(98,244)(63,244)(50,269)
\qbezier(50,269)(84,271)(98,244)
\put(50,286){\circle*{4}}
\qbezier(169,299)(185,324)(217,324)
\put(147,297){$v_{1}$}
\put(150,354){$v_{2}$}
\put(91,372){$v_{3}$}
\put(87,234){$v_{g-1}$}
\put(146,247){$v_{g}$}
\put(470,272){\circle*{4}}
\put(470,338){\circle*{4}}
\qbezier(470,272)(452,303)(470,338)
\qbezier(470,338)(496,307)(470,272)
\put(522,366){\circle*{4}}
\qbezier(470,338)(487,367)(522,366)
\qbezier(470,338)(509,335)(522,366)
\put(579,341){\circle*{4}}
\qbezier(522,366)(568,368)(579,341)
\qbezier(522,366)(541,336)(579,341)
\put(505,243){\circle*{4}}
\qbezier(470,272)(475,248)(505,243)
\qbezier(470,272)(501,274)(505,243)
\put(553,246){\circle*{4}}
\qbezier(505,243)(533,235)(553,246)
\qbezier(505,243)(530,267)(553,246)
\put(581,284){\circle*{4}}
\qbezier(553,246)(583,257)(581,284)
\qbezier(553,246)(555,278)(581,284)
\put(468,306){\circle*{4}}
\qbezier(468,306)(443,332)(419,326)
\put(579,300){\circle*{4}}
\put(579,312){\circle*{4}}
\put(579,324){\circle*{4}}
\put(452,270){$v_{1}$}
\put(452,340){$v_{2}$}
\put(516,372){$v_{3}$}
\put(498,233){$v_{g}$}
\put(547,234){$v_{g-1}$}
\put(473,305){$v_{a(1,1)}$}
\put(197,297){$T_{1}$}
\put(424,300){$T_{2}$}
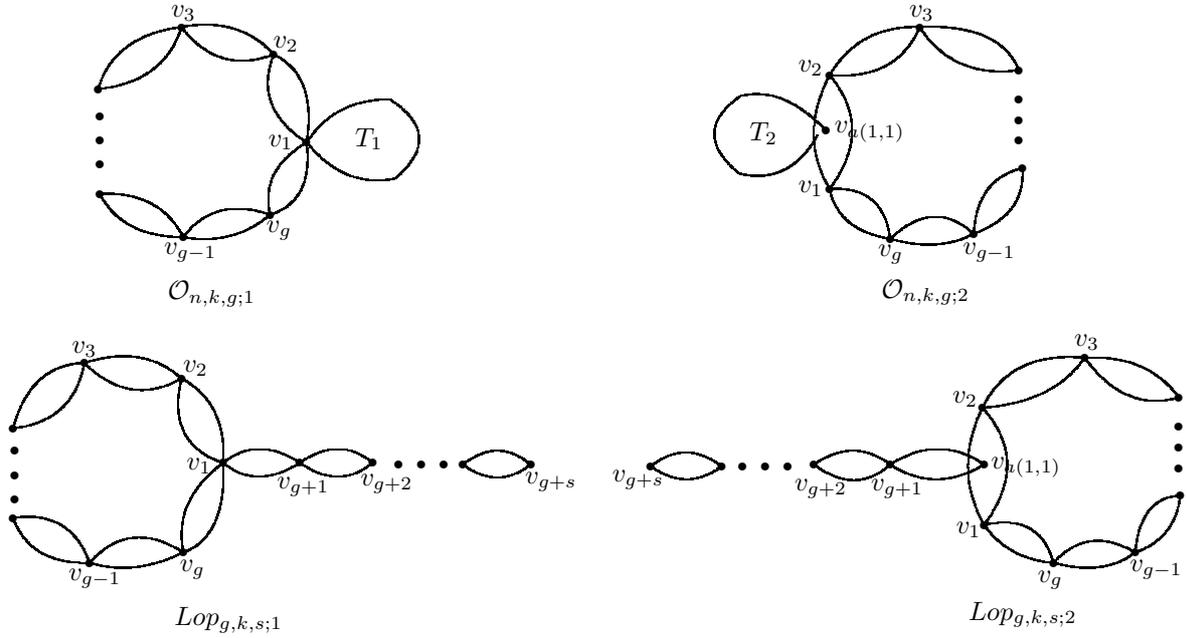
\put(190,-9){Fig. 1.1. $\mathcal{O}_{n,k,g;1}$, $\mathcal{O}_{n,k,g;2}$, $Lop_{g,k,s;1}$, $Lop_{g,k,s;2}$}
\qbezier(169,299)(191,272)(220,278)
\qbezier(217,324)(249,303)(220,278)
\put(89,208){$\mathcal{O}_{n,k,g;1}$}
\qbezier(463,303)(447,274)(418,281)
\qbezier(419,326)(389,305)(418,281)
\put(500,209){$\mathcal{O}_{n,k,g;2}$}
\put(97,162){\circle*{4}}
\put(121,113){\circle*{4}}
\qbezier(97,162)(123,148)(121,113)
\put(98,61){\circle*{4}}
\qbezier(121,113)(123,70)(98,61)
\put(41,171){\circle*{4}}
\qbezier(97,162)(75,182)(41,171)
\put(0,133){\circle*{4}}
\qbezier(41,171)(10,170)(0,133)
\put(44,55){\circle*{4}}
\qbezier(98,61)(70,47)(44,55)
\put(0,81){\circle*{4}}
\qbezier(44,55)(10,54)(0,81)
\qbezier(41,171)(35,137)(0,133)
\qbezier(41,171)(67,145)(97,162)
\qbezier(97,162)(88,124)(121,113)
\qbezier(121,113)(92,91)(98,61)
\qbezier(98,61)(67,81)(44,55)
\qbezier(0,81)(33,86)(44,55)
\put(1,120){\circle*{4}}
\put(1,106){\circle*{4}}
\put(1,92){\circle*{4}}
\put(165,113){\circle*{4}}
\qbezier(121,113)(141,129)(165,113)
\qbezier(121,113)(140,97)(165,113)
\put(207,113){\circle*{4}}
\qbezier(165,113)(188,129)(207,113)
\qbezier(165,113)(188,97)(207,113)
\put(248,112){\circle*{4}}
\put(235,112){\circle*{4}}
\put(222,112){\circle*{4}}
\put(259,112){\circle*{4}}
\put(298,112){\circle*{4}}
\qbezier(259,112)(279,128)(298,112)
\qbezier(259,112)(279,97)(298,112)
\put(558,145){\circle*{4}}
\put(559,77){\circle*{4}}
\qbezier(558,145)(541,109)(559,77)
\put(617,174){\circle*{4}}
\qbezier(558,145)(574,178)(617,174)
\put(671,151){\circle*{4}}
\qbezier(617,174)(653,176)(671,151)
\put(599,55){\circle*{4}}
\qbezier(559,77)(573,57)(599,55)
\put(646,61){\circle*{4}}
\qbezier(599,55)(628,47)(646,61)
\put(672,94){\circle*{4}}
\qbezier(646,61)(669,69)(672,94)
\qbezier(558,145)(586,114)(559,77)
\qbezier(558,145)(599,148)(617,174)
\qbezier(617,174)(640,142)(671,151)
\qbezier(559,77)(591,84)(599,55)
\qbezier(599,55)(619,77)(646,61)
\qbezier(672,94)(649,93)(646,61)
\put(671,110){\circle*{4}}
\put(671,122){\circle*{4}}
\put(671,134){\circle*{4}}
\put(505,112){\circle*{4}}
\put(559,112){\circle*{4}}
\qbezier(505,112)(533,130)(559,112)
\qbezier(505,112)(531,95)(559,112)
\put(461,112){\circle*{4}}
\qbezier(505,112)(484,128)(461,112)
\qbezier(461,112)(485,94)(505,112)
\put(367,111){\circle*{4}}
\put(408,111){\circle*{4}}
\put(446,111){\circle*{4}}
\put(433,111){\circle*{4}}
\put(420,111){\circle*{4}}
\qbezier(367,111)(387,127)(408,111)
\qbezier(367,111)(388,96)(408,111)
\put(93,19){$Lop_{g,k,s;1}$}
\put(550,22){$Lop_{g,k,s;2}$}
\put(100,109){$v_{1}$}
\put(98,166){$v_{2}$}
\put(34,177){$v_{3}$}
\put(96,50){$v_{g}$}
\put(32,45){$v_{g-1}$}
\put(152,98){$v_{g+1}$}
\put(200,100){$v_{g+2}$}
\put(295,101){$v_{g+s}$}
\put(543,72){$v_{1}$}
\put(563,110){$v_{a(1,1)}$}
\put(541,148){$v_{2}$}
\put(611,181){$v_{3}$}
\put(590,44){$v_{g}$}
\put(643,50){$v_{g-1}$}
\put(494,97){$v_{g+1}$}
\put(450,97){$v_{g+2}$}
\put(345,102){$v_{g+s}$}
\end{picture}
\end{center}

Let $E_{ij} = \{e \in E : i, j \in e\}$. The $adjacency$ $matrix$ $\mathcal {A}_{\mathcal{G}}=[(\mathcal {A}_{\mathcal{G}})_{ij}]$ of a hypergraph $\mathcal{G}$ is defined as
$$(\mathcal {A}_{\mathcal{G}})_{ij}=\sum_{e\in E_{ij}}\frac{1}{|e| - 1}.$$ It is easy to find that $\mathcal {A}_{\mathcal{G}}$ is symmetric if there is no requirement for direction on hypergraph $\mathcal{G}$, and find that $\mathcal {A}_{\mathcal{G}}$ is very convenient to be used to investigate the spectum of a hypergraph even without the requirement for edge uniformity. The $spectral$ $radius$ $\rho(\mathcal {G})$ of a hypergraph $\mathcal{G}$ is defined to be the spectral radius $\rho(\mathcal {A}_{\mathcal{G}})$, which is the maximum modulus among all eigenvalues of $\mathcal {A}_{\mathcal{G}}$.
In spectral theory of hypergraphs, the spectral radius is an index that attracts much attention due to its fine properties \cite{KCKPZ, YFTP, SFSH, LSQ, HLBZ, LLSM, COQY, YQY}.

We assume that the hypergraphs throughout
this paper are simple, i.e. $e_{i} \neq e_{j}$ if $i \neq j$, and assume the hypergraphs throughout
this paper are undirected.

Let $G = (V, E)$ be an ordinary graph (2-uniform). For every $k \geq 3$,
the $k$th power of $G$, denoted by $G^{k}= (V^{k} , E^{k} )$, is defined as the $k$-uniform hypergraph with
the edge set $E^{k}= \{e \cup \{v_{e_{1}}$, $v_{e_{2}}$, $\ldots$, $v_{e_{k-2}}\}:e\in E\}$ and the vertex set $V^{k}= V \cup
(\cup_{e\in E} \{v_{e_{1}}$, $v_{e_{2}}$, $\ldots$, $v_{e_{k-2}}\})$, where $V \cap
(\cup_{e\in E} \{v_{e_{1}}$, $v_{e_{2}}$, $\ldots$, $v_{e_{k-2}}\})=\emptyset$, $\{v_{e_{1}}$, $v_{e_{2}}$, $\ldots$, $v_{e_{k-2}}\}\cap \{v_{f_{1}}$, $v_{f_{2}}$, $\ldots$, $v_{f_{k-2}}\}=\emptyset$ for $e\neq f$, $e,f\in E$. For an ordinary graph $G$, we call it regular or $d$-regular if its very vertex degree is constant $d$. The following Theorem \ref{le3,28} establish a relation between the sepectral radius of an ordinary graph and the sepectral radius of a hypergraph.

\begin{theorem} \label{le3,28}
Suppose $\mathcal{G}=G^{k}$ is the kth power of degree $d$-regular graph $G$ ($k\geq 3$). Then $$\rho(\mathcal{G})=\frac{k-3+d+\sqrt{(k-3+d)^{2}+4d(k-1)}}{2(k-1)}.$$
\end{theorem}

As well, in this paper, for $k$-uniform ($k\geq 3$) unicyclic hypergraphs with fixed number of vertices, for $k$-uniform ($k\geq 3$) unicyclic hypergraphs with fixed number of vertices and fixed girth, for $k$-uniform ($k\geq 3$) octopuslike hypergraphs, and for $k$-uniform ($k\geq 3$) $lollipop$ hypergraphs, we get the following result:

\begin{theorem}\label{th01.03} %------
(1) Let $\mathcal{G}$ be a $k$-uniform ($k\geq 3$) unicyclic hypergraph of order $n$. Then $$\frac{k-1+\sqrt{(k-1)^{2}+8(k-1)}}{2(k-1)}\leq \rho(\mathcal{G})$$ with equality if and only if $\mathcal{G}\cong \mathscr{C}(n, k)$.

(2) Let $\mathcal{G}$ be a $k$-uniform ($k\geq 3$) unicyclic hypergraph of order $n$ satisfying that $\mathcal{G}\ncong \mathscr{C}(n, k)$. Then $\rho(Lop_{g,k,1;2})\leq \rho(\mathcal{G})$ with equality if and only if $\mathcal{G}\cong Lop_{g,k,1;2}$, where $g=\frac{n-k+1}{k-1}$.

\end{theorem}

\begin{theorem}\label{th01.01} %------
Let $\mathcal{G}$ be a $k$-uniform ($k\geq 3$) unicyclic hypergraph of order $n$ and girth $g$. Then $\rho(\mathcal{G})\leq\rho(\mathscr{U}^{\ast}(n,k,g))$ with equality if and only if $\mathcal{G}\cong \mathscr{U}^{\ast}(n,k,g)$.

\end{theorem}

\begin{theorem}\label{th01.06} %------
(1) Let $\mathcal{G}$ be a $k$-uniform ($k\geq 3$) octopuslike hypergraph $\mathcal{O}_{n,k,g}$ of order $n$. Then $\rho(Lop_{g,k,s;2})\leq \rho(\mathcal{G})$ with equality if and only if $\mathcal{G}\cong Lop_{g,k,s;2}$ where $s=\frac{n}{k-1}-g$.

(2) Let $\mathcal{G}$ be a $k$-uniform ($k\geq 3$) octopuslike hypergraph $\mathcal{O}_{n,k,g;1}$ of order $n$. Then $\rho(Lop_{g,k,s;1})\leq \rho(\mathcal{G})$ with equality if and only if $\mathcal{G}\cong Lop_{g,k,s;1}$ where $s=\frac{n}{k-1}-g$.

(3) Let $\mathcal{G}$ be a $k$-uniform ($k\geq 3$) octopuslike hypergraph $\mathcal{O}_{n,k,g;2}$ of order $n$. Then $\rho(Lop_{g,k,s;2})\leq \rho(\mathcal{G})$ with equality if and only if $\mathcal{G}\cong Lop_{g,k,s;2}$ where $s=\frac{n}{k-1}-g$.

\end{theorem}

\begin{theorem}\label{th01.07} %------
(1) $\rho(Lop_{g+1,k,s;1})<\rho(Lop_{g,k,s;1})$ where $s\geq 1$.

(2) $\rho(Lop_{g+1,k,s;2})<\rho(Lop_{g,k,s;2})$ where $s\geq 1$.

\end{theorem}

The layout of this paper is as follows: section 2 introduces some basic knowledge; section 3 represents our results.

\section{Preliminary}

\ \ \ \ \ For the requirements in the narrations afterward, we need some prepares. For a hypergraph $\mathcal{G}$ with vertex set $\{v_{1}$, $v_{2}$, $\ldots$, $v_{n}\}$, a vector $X=(x_{v_1}, x_{v_2}, \ldots, x_{v_n})^T \in R^n$ on $\mathcal{G}$ is a vector that
 entry $x_{v_i}$ is mapped to vertex $v_i$ for $i\leq i\leq n$.

 From \cite{OP}, by the famous Perron-Frobenius theorem, for $\mathcal {A}_{\mathcal{G}}$ of a connected uniform hypergraph $\mathcal{G}$ of order $n$, we know that there is unique one positive eigenvector $X=(x_{v_{1}}$, $x_{v_{2}}$, $\ldots$, $x_{v_{n}})^T \in R^{n}_{++}$ ($R^{n}_{++}$ means the set of positive real vectors of dimension $n$) corresponding to $\rho(\mathcal{G})$, where $\sum^{n}_{i=1}x^{2}_{v_{i}}= 1$ and each entry $x_{v_i}$ is mapped to each vertex $v_i$ for $i\leq i\leq n$. We call such an eigenvector $X$
the $principal$ $eigenvector$ of $\mathcal{G}$.

Let $A$ be an irreducible nonnegative $n \times n$ real matrix (with every entry being real number) with spectral radius $\rho(A)$. The following extremal representation (Rayleigh quotient) will be useful:
$$\rho(A)=\max_{X\in R^{n}, X\neq0}\frac{X^{T}AX}{X^{T}X},$$ and if a vector $X$ satisfies that $\frac{X^{T}AX}{X^{T}X}=\rho(A)$, then $AX=\rho(A)X$.

\section{Main results}

\begin{lemma}{\bf \cite{{GYSLHZ}}}\label{le03,00} %------
Let $G$ be a $k$-uniform hypergraph on $n$ vertices. A permutation $\sigma \in S_{n}$ (sysmmetric group of degree $n$) is an automorphism of
$G$ if and only if $P_{\sigma}\mathcal {A}_{G} =\mathcal {A}_{G}P_{\sigma}$.
\end{lemma}

\begin{lemma}{\bf \cite{{GYSLHZ}}}\label{le03,01} %------
Let $G$ be a connected $k$-uniform hypergraph, $\mathcal {A}_{G}$
be its
adjacency matrix. If $\mathcal {A}X=\lambda X$, then for each automorphism $\sigma$ of $G$, we have $\mathcal {A} (P_{\sigma} X)=\lambda (P_{\sigma}X)$.
Moreover, if $\mathcal {A}X=\rho(\mathcal {A}) X$, then

$\mathrm{(1)}$ $P_{\sigma} X = X$;

$\mathrm{(2)}$ for any orbit $\Omega$ of $Aut (G)$ and each pair of vertices $i, j\in \Omega$, $x_{i}=x_{j}$ in $X$.
\end{lemma}

\begin{lemma}{\bf \cite{GYSLZ}}\label{le03.02} %------
Let $e_{1}=\{v_{1,1}$, $v_{1,2}$, $\ldots$, $v_{1,k}\}$, $e_{2}=\{v_{2,1}$, $v_{2,2}$, $\ldots$, $v_{2,k}\}$, $\ldots$, $e_{j}=\{v_{j,1}$, $v_{j,2}$, $\ldots$, $v_{j,k}\}$ be some edges in a connected $k$-uniform hypergraph $\mathcal{G}$; $v_{u,1}$, $v_{u,2}$, $\ldots$, $v_{u,t}$ be vertices in $\mathcal{G}$ that $t< k$. For $1\leq i\leq j$, $\{v_{u,1}, v_{u,2}, \ldots, v_{u,t}\}\nsubseteq e_{i}$, $e^{'}_{i}=(e_{i}\setminus \{v_{i,1}, v_{i,2}, \ldots, v_{i,t}\})\cup\{v_{u,1}, v_{u,2}, \ldots, v_{u,t}\}$ satisfying that $e^{'}_{i}\notin E(\mathcal{G})$. Let $\mathcal{G}^{'}=\mathcal{G}-\sum e_{i}+\sum e^{'}_{i}$. If in the principal eigenvector $X$ of $\mathcal{G}$, for $1\leq i\leq j$, $x_{v_{i,1}}\leq x_{v_{u,1}}$, $x_{v_{i,2}}\leq x_{v_{u,2}}$, $\ldots$, and $x_{v_{i,t}}\leq x_{v_{u,t}}$, then $\rho(\mathcal{G}^{'})>\rho(\mathcal{G})$.
\end{lemma}

\begin{lemma}{\bf \cite{GYSLZ}} \label{le3,04,01}
Let $A$ be an irreducible nonnegative square symmetric real matrix with order $n$
and spectral radius $\rho$, $y\in R^{n}_{++}$ be a positive vector. If there exists $r\in R_{+}$ such that $Ay\leq ry$, then
$\rho\leq r$ with equality if and only if $Ay= ry$. Similarly, if there exists $r\in R_{+}$ such that $Ay\geq ry$, then
$\rho\geq r$ with equality if and only if $Ay= ry$.
\end{lemma}

\

\setlength{\unitlength}{0.6pt}
\begin{center}
\begin{picture}(578,116)
\qbezier(0,60)(0,72)(28,81)\qbezier(28,81)(57,90)(98,90)\qbezier(98,90)(138,90)(167,81)\qbezier(167,81)(196,72)(196,60)\qbezier(196,60)(196,47)(167,38)
\qbezier(167,38)(138,30)(98,30)\qbezier(98,30)(57,30)(28,38)\qbezier(28,38)(0,47)(0,60)
\qbezier(132,60)(132,74)(176,85)\qbezier(176,85)(221,95)(285,95)\qbezier(285,95)(348,95)(393,85)\qbezier(393,85)(438,74)(438,60)\qbezier(438,60)(438,45)(393,34)
\qbezier(393,34)(348,24)(285,24)\qbezier(285,25)(221,25)(176,34)\qbezier(176,34)(131,45)(132,60)
\qbezier(373,60)(373,72)(402,80)\qbezier(402,80)(432,89)(475,89)\qbezier(475,89)(517,89)(547,80)\qbezier(547,80)(577,72)(577,60)\qbezier(577,60)(577,47)(547,39)
\qbezier(547,39)(517,30)(475,30)\qbezier(475,31)(432,31)(402,39)\qbezier(402,39)(372,47)(373,60)
\put(310,65){\circle*{4}}
\put(299,65){\circle*{4}}
\put(164,65){\circle*{4}}
\put(157,51){$v_{1}$}
\put(289,65){\circle*{4}}
\put(252,65){\circle*{4}}
\put(221,65){\circle*{4}}
\put(406,66){\circle*{4}}
\put(399,51){$v_{k}$}
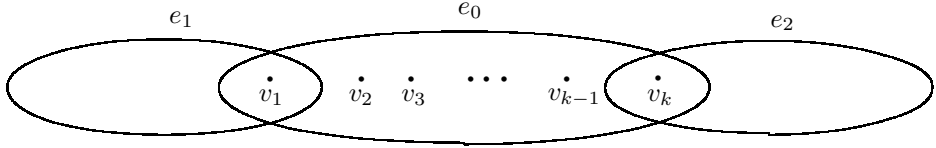
\put(195,-9){Fig. 3.1. $e_{0}$, $e_{1}$, $e_{2}$ in $\mathcal{G}$}
\put(101,100){$e_{1}$}
\put(281,106){$e_{0}$}
\put(476,98){$e_{2}$}
\put(213,51){$v_{2}$}
\put(246,51){$v_{3}$}
\put(349,65){\circle*{4}}
\put(337,52){$v_{k-1}$}
\end{picture}
\end{center}

\begin{PProof}
For any edge $e_{i}=uv\in E(G)$, we denote by $e^{'}_{i}=\{u, w_{i_{1}}, w_{i_{2}}, \ldots, w_{i_{k-2}}, v\}$ the corresponding edge in $\mathcal{G}$. Let $t=\frac{k-3+d+\sqrt{(k-3+d)^{2}+4d(k-1)}}{2(k-1)}$, $Y$ be a vector on $\mathcal{G}$ satisfying that
$$\left \{\begin{array}{ll}
y_{s}=\frac{2}{(k-1)t-(k-3)}\ & \ s\in\{w_{i_{1}}, w_{i_{2}}, \ldots, w_{i_{k-2}}\}, e^{'}_{i}\in E(\mathcal{G})\\
\\
 y_{s}=1,\ & \ others.\end{array}\right.$$

Then for any vertex $u\in V(G)$, we have $$(\mathcal {A}_{\mathcal{G}}Y)_{u}=\sum_{uv\in E(G)}\frac{\frac{2(k-2)}{(k-1)t-(k-3)}+y_{v}}{(k-1)}=d\cdot\frac{\frac{2(k-2)}{(k-1)t-(k-3)}+1}{(k-1)}=ty_{u}=t.$$
For any vertex $z\in\{w_{i_{1}}, w_{i_{2}}, \ldots, w_{i_{k-2}}\}$ where $e^{'}_{i}=\{u, w_{i_{1}}, w_{i_{2}}, \ldots, w_{i_{k-2}}, v\}\in E(\mathcal{G})$, we have $$(\mathcal {A}_{\mathcal{G}}Y)_{z}=\frac{\frac{2(k-3)}{(k-1)t-(k-3)}+x_{u}+x_{v}}{(k-1)}=\frac{\frac{2(k-3)}{(k-1)t-(k-3)}+2}{(k-1)}=ty_{z}.$$ Thus it follows that
$\mathcal {A}_{\mathcal{G}}Y=tY$. Then the result follows. This completes the proof. \ \ \ \ \ $\Box$
\end{PProof}

\begin{corollary} \label{cl01,04}
$\rho(\mathscr{C}(n, k))=\frac{k-1+\sqrt{(k-1)^{2}+8(k-1)}}{2(k-1)}$.
\end{corollary}

\begin{proof}
Note that the ordinary cycle is 2-regular graph, and $\mathscr{C}(n, k)$ is the kth power of a ordinary cycle with length $\frac{n}{k-1}$. Then this corollary follows from Theorem \ref{le3,28}. This completes the proof. \ \ \ \ \ $\Box$
\end{proof}

\setlength{\unitlength}{0.7pt}
\begin{center}
\begin{picture}(632,214)
\put(270,-9){Fig. 3.2. $\mathcal{G}_{1}$, $\mathcal{G}_{2}$}
\put(483,170){$v_{1}$}
\put(456,181){$e_{1}$}
\put(430,170){$v_{2}$}
\put(398,181){$e_{4}$}
\put(508,182){$e_{2}$}
\put(551,180){$e_{3}$}
\put(457,49){$v_{1}$}
\put(427,61){$e_{1}$}
\put(405,49){$v_{2}$}
\put(379,62){$e_{4}$}
\put(485,58){$e'$}
\put(515,50){$u$}
\put(538,58){$e'_{2}$}
\put(561,47){$v_{3}$}
\put(584,60){$e_{3}$}
\put(343,62){\circle*{4}}
\put(354,62){\circle*{4}}
\put(365,62){\circle*{4}}
\put(581,181){\circle*{4}}
\put(592,181){\circle*{4}}
\put(603,181){\circle*{4}}
\put(367,181){\circle*{4}}
\put(378,181){\circle*{4}}
\put(389,181){\circle*{4}}
\put(610,61){\circle*{4}}
\put(621,61){\circle*{4}}
\put(632,61){\circle*{4}}
\qbezier(485,151)(485,128)(485,105)
\qbezier(497,151)(497,128)(497,105)
\qbezier(476,117)(484,106)(492,95)
\qbezier(492,95)(498,107)(505,118)
\put(479,15){$\mathcal{G}_{2}$}
\qbezier(580,6)(581,6)(580,6)
\qbezier(142,148)(142,125)(142,101)
\qbezier(153,148)(153,125)(153,101)
\qbezier(134,114)(141,103)(148,91)
\qbezier(162,114)(155,103)(148,91)
\put(266,58){\circle*{4}}
\put(277,58){\circle*{4}}
\put(288,58){\circle*{4}}
\put(57,185){$e_{4}$}
\put(86,177){$v_{2}$}
\put(115,186){$e_{1}$}
\put(146,175){$v_{1}$}
\put(174,185){$e_{2}$}
\put(202,175){$v_{3}$}
\put(225,184){$e_{3}$}
\put(34,59){$e_{4}$}
\put(60,49){$v_{2}$}
\put(87,58){$e'_{1}$}
\put(117,52){$u$}
\put(144,56){$e'$}
\put(167,50){$v_{1}$}
\put(193,59){$e_{2}$}
\put(215,49){$v_{3}$}
\put(0,59){\circle*{4}}
\put(11,59){\circle*{4}}
\put(22,59){\circle*{4}}
\put(23,184){\circle*{4}}
\put(34,184){\circle*{4}}
\put(45,184){\circle*{4}}
\put(238,58){$e_{3}$}
\put(250,184){\circle*{4}}
\put(261,183){\circle*{4}}
\put(272,183){\circle*{4}}
\put(135,16){$\mathcal{G}_{1}$}
\put(89,189){\circle*{4}}
\put(209,186){\circle*{4}}
\put(487,181){\circle*{4}}
\put(432,181){\circle*{4}}
\put(225,63){\circle*{4}}
\put(171,62){\circle*{4}}
\put(120,62){\circle*{4}}
\put(62,63){\circle*{4}}
\put(517,60){\circle*{4}}
\put(409,62){\circle*{4}}
\put(150,187){\circle*{4}}
\qbezier(89,189)(122,161)(150,187)
\qbezier(89,189)(124,214)(150,187)
\qbezier(150,187)(183,212)(209,186)
\qbezier(150,187)(183,163)(209,186)
\qbezier(209,186)(223,202)(244,200)
\qbezier(209,186)(225,165)(244,169)
\qbezier(56,200)(78,202)(89,189)
\qbezier(89,189)(77,170)(54,171)
\qbezier(62,63)(95,85)(120,62)
\qbezier(62,63)(92,38)(120,62)
\qbezier(120,62)(151,84)(171,62)
\qbezier(120,62)(150,39)(171,62)
\qbezier(171,62)(200,82)(225,63)
\qbezier(171,62)(201,39)(225,63)
\qbezier(225,63)(229,77)(256,75)
\qbezier(225,63)(226,40)(257,44)
\qbezier(62,63)(58,76)(32,74)
\qbezier(62,63)(58,42)(32,46)
\qbezier(432,181)(466,207)(487,181)
\put(569,60){\circle*{4}}
\put(538,182){\circle*{4}}
\qbezier(432,181)(467,158)(487,181)
\qbezier(487,181)(516,207)(538,182)
\qbezier(487,181)(516,160)(538,182)
\qbezier(538,182)(544,199)(569,198)
\qbezier(538,182)(547,162)(571,166)
\qbezier(432,181)(421,199)(400,198)
\qbezier(432,181)(418,162)(398,166)
\put(462,62){\circle*{4}}
\qbezier(409,62)(437,85)(462,62)
\qbezier(462,62)(491,86)(517,60)
\qbezier(517,60)(549,85)(569,60)
\qbezier(409,62)(434,42)(462,62)
\qbezier(462,62)(490,43)(517,60)
\qbezier(517,60)(549,42)(569,60)
\qbezier(409,62)(405,77)(379,78)
\qbezier(409,62)(398,48)(379,49)
\qbezier(569,60)(571,73)(600,77)
\qbezier(569,60)(573,45)(603,47)
\put(531,170){$v_{3}$}
\end{picture}
\end{center}

\begin{lemma} {\bf \cite{GYSLZ}}\label{le3,13}
Suppose $\mathcal{G}$ is a connected hypergraph with spectral radius $\rho$ and principal eigenvector $X$. $e_{1}$, $e_{2}$, $e_{3}$, $e_{4}$ are edges in $\mathcal{G}$, where $|e_{1}|, |e_{2}|, |e_{3}|, |e_{4}|\geq 3$, $e_{1}\cap e_{2}=\{v_{1}\}$, $e_{1}\cap e_{4}=\{v_{2}\}$, $e_{2}\cap e_{3}=\{v_{3}\}$, $deg_{\mathcal{G}}(v_{1})=deg_{\mathcal{G}}(v_{2})=deg_{\mathcal{G}}(v_{3})=2$,
$deg_{\mathcal{G}}(v)=1$ for $v\in (e_{1}\cup e_{2})\setminus \{v_{1}, v_{2}, v_{3}\}$.
Let $e^{'}_{1}=(e_{1}\setminus \{v_{1}\})\cup \{u\}$, $e^{'}_{2}=(e_{2}\setminus \{v_{1}\})\cup \{u\}$, $e^{'}=\{v_{1}$, $u_{1}$, $u_{2}$, $\ldots$, $u_{t-2}$, $u\}$ where $u\notin V(\mathcal{G})$, $u_{i}\notin V(\mathcal{G})$ for $1\leq i\leq t-2$, $\mathcal{G}_{1}=\mathcal{G}-e_{1}+e^{'}_{1}+e^{'}$, $\mathcal{G}_{2}=\mathcal{G}-e_{2}+e^{'}_{2}+e^{'}$ (see Fig. 3.2).
If $t\geq \max\{e_{1}$, $e_{2}\}$, $x_{v_{1}}\leq x_{v_{2}}$, $x_{v_{1}}\leq x_{v_{3}}$, then $\rho(\mathcal{G}_{1})\leq\rho(\mathcal{G})$, $\rho(\mathcal{G}_{2})\leq\rho(\mathcal{G})$ with either equality holding if and only if $|e^{'}|=|e_{1}|= |e_{2}|$, $x_{v_{1}}=x_{v_{2}}=x_{v_{3}}=x_{u}$ and $x_{z}=x_{\omega}$ for $z,\omega\in (e_{1}\cup e_{2}\cup e^{'})\setminus\{v_{1}, v_{2}, u\}$.

\end{lemma}

\begin{lemma}{\bf \cite{GYSLZ}} \label{le3,14}
Let $e$ be a new edge not containing in connected hypergraph $\mathcal{G}$. Let $\mathcal{G}^{'}=\mathcal{G}+e$. If $\mathcal{G}^{'}$ is also connected, then $\rho(\mathcal{G}^{'})>\rho(\mathcal{G})$.
\end{lemma}

\setlength{\unitlength}{0.7pt}
\begin{center}
\begin{picture}(461,211)
\put(1,140){\circle*{4}}
\put(1,126){\circle*{4}}
\put(0,162){\circle*{4}}
\put(48,192){\circle*{4}}
\qbezier(0,162)(8,188)(48,192)
\qbezier(0,162)(38,160)(48,192)
\put(101,176){\circle*{4}}
\qbezier(48,192)(82,198)(101,176)
\qbezier(48,192)(69,164)(101,176)
\put(120,125){\circle*{4}}
\qbezier(101,176)(126,164)(120,125)
\qbezier(101,176)(89,145)(120,125)
\put(99,80){\circle*{4}}
\qbezier(120,125)(124,91)(99,80)
\qbezier(120,125)(93,109)(99,80)
\put(48,65){\circle*{4}}
\qbezier(99,80)(79,61)(48,65)
\qbezier(48,65)(67,94)(99,80)
\put(1,90){\circle*{4}}
\qbezier(48,65)(14,65)(1,90)
\qbezier(1,90)(35,97)(48,65)
\put(1,112){\circle*{4}}
\qbezier(120,125)(136,150)(168,150)
\put(101,124){$v_{1}$}
\put(101,179){$v_{2}$}
\put(42,198){$v_{3}$}
\put(38,55){$v_{g-1}$}
\put(98,70){$v_{g}$}
\put(351,97){\circle*{4}}
\put(350,163){\circle*{4}}
\qbezier(351,97)(332,128)(350,163)
\qbezier(350,163)(376,132)(351,97)
\put(400,191){\circle*{4}}
\qbezier(350,163)(367,192)(400,191)
\qbezier(350,163)(389,160)(400,191)
\put(458,173){\circle*{4}}
\qbezier(400,191)(434,196)(458,173)
\qbezier(400,191)(421,161)(458,173)
\put(385,68){\circle*{4}}
\qbezier(351,97)(357,77)(385,68)
\qbezier(351,97)(381,99)(385,68)
\put(433,69){\circle*{4}}
\qbezier(385,68)(413,60)(433,69)
\qbezier(385,68)(410,92)(433,69)
\put(461,103){\circle*{4}}
\qbezier(433,69)(460,80)(461,103)
\qbezier(433,69)(435,103)(461,103)
\put(348,131){\circle*{4}}
\qbezier(348,131)(323,157)(299,151)
\put(459,127){\circle*{4}}
\put(459,139){\circle*{4}}
\put(459,151){\circle*{4}}
\put(336,93){$v_{1}$}
\put(333,163){$v_{2}$}
\put(392,198){$v_{3}$}
\put(379,57){$v_{g}$}
\put(429,58){$v_{g-1}$}
\put(352,129){$v_{a(1,1)}$}
\put(148,123){$G_{1}$}
\put(304,125){$G_{2}$}
\put(180,-9){Fig. 3.3. $\mathbb{G}_{1,g}, \mathbb{G}_{2,g}$}
\qbezier(120,125)(142,98)(171,104)
\qbezier(168,150)(200,129)(171,104)
\put(46,27){$\mathbb{G}_{1,g}$}
\qbezier(343,128)(327,99)(298,106)
\qbezier(299,151)(269,130)(298,106)
\put(388,28){$\mathbb{G}_{2,g}$}
\end{picture}
\end{center}

Let $G$ be a $k$-uniform hypergraph which is disjoint to $\mathcal{C}_{g}$ ($\mathcal{C}_{g}$ is depicted in section 1), that is, they have no common vertices. Let $\mathbb{G}_{1,g}$ be the $k$-uniform hypergraph obtained by identifying the vertex $u$ of $G$ and the vertex $v_{1}$ of $\mathcal{C}_{g}$, where we still denote by $v_{1}$ the new vertex obtained by identifying $u$ and $v_{1}$ (see Fig. 3.3), denote by $G_{1}$ obtained from $G$ after identifying $u$ and $v_{1}$; let $\mathbb{G}_{2,g}$ be the $k$-uniform hypergraph obtained by identifying the vertex $u$ of $G$ and the vertex $v_{a(1,1)}$ of $\mathcal{C}_{g}$, where we still denote by $v_{a(1,1)}$ the new vertex obtained by identifying $u$ and $v_{a(1,1)}$ (see Fig. 3.3), denote by $G_{2}$ obtained from $G$ after identifying $u$ and $v_{a(1,1)}$.

\begin{lemma} \label{le3,19}
(1) Let $X$ be
the principal eigenvector of $\mathbb{G}_{1,g}$.

(1.1) If $g$ is odd, then $x_{v_{1+i}}=x_{v_{g-i+1}}$, $x_{v_{a(i,1)}}=x_{v_{a(g-i+1,1)}}$, $x_{v_{i}}\geq x_{v_{i+1}}$ for $i=1, 2, \ldots, \frac{g-1}{2}$.

(1.2) If $g$ is even, then $x_{v_{1+i}}=x_{v_{g-i+1}}$, $x_{v_{a(i,1)}}=x_{v_{a(g-i+1,1)}}$, $x_{v_{i}}\geq x_{v_{i+1}}$ for $i=1, 2, \ldots, \frac{g-2}{2}$.

(2) Let $X$ be
the principal eigenvector of $\mathbb{G}_{2,g}$.

(2.1) If $g$ is odd, then $x_{v_{1}}=x_{v_{2}}$, $x_{v_{2+i}}=x_{v_{g-i+1}}$ for $i=1, 2, \ldots, \frac{g-3}{2}$, $x_{v_{a(1+i,1)}}=x_{v_{a(g-i+1,1)}}$ for $i=1, 2, \ldots, \frac{g-1}{2}$, $x_{v_{i}}\geq x_{v_{i+1}}$ for $i=2, 3, \ldots, \frac{g+1}{2}$.

(2.2) If $g$ is even, then $x_{v_{1}}=x_{v_{2}}$, $x_{v_{2+i}}=x_{v_{g-i+1}}$ for $i=1, 2, \ldots, \frac{g-2}{2}$, $x_{v_{a(1+i,1)}}=x_{v_{a(g-i+1,1)}}$ for $i=1, 2, \ldots, \frac{g-2}{2}$, $x_{v_{i}}\geq x_{v_{i+1}}$ for $i=2, 3, \ldots, \frac{g}{2}$.
\end{lemma}

\begin{proof}
(1) By Lemma \ref{le03,01}, we know that $x_{v_{a(i,j)}}=x_{v_{a(i,1)}}$ for $i=1, 2, \ldots, g$, $j=1, 2, \ldots, k-2$.

(1.1) Using Lemma \ref{le03,01} gets $x_{v_{1+i}}=x_{v_{g-i+1}}$, $x_{v_{a(i,1)}}=x_{v_{a(g-i+1,1)}}$ for $i=1, 2, \ldots, \frac{g-1}{2}$.
Next, we prove $x_{v_{i}}\geq x_{v_{i+1}}$ for $i=1, 2, \ldots, \frac{g-1}{2}$. We prove it by contradiction. Suppose that there exists some $2\leq i \leq\frac{g+1}{2}$ such that $x_{v_{i}}> x_{v_{i-1}}$. Along the hyperpath $P_{1}=v_{\frac{g+1}{2}}e_{\frac{g-1}{2}}v_{\frac{g-1}{2}}\cdots e_{2}v_{2}e_{1}v_{1}$, suppose $\eta$ is the largest number form $2$ to $\frac{g+1}{2}$ such that $x_{v_{\eta}}> x_{v_{\eta-1}}$.

{\bf Case 1} There is some number $2\leq\tau<\eta$ that $x_{v_{\eta}}> x_{v_{i}}$ for $\tau\leq i\leq\eta-1$, but $x_{v_{\eta}}\leq x_{v_{\tau-1}}$.

Let $Y$ be
the vector on $\mathbb{G}_{1,g}$ satisfying that
$$\left \{\begin{array}{ll}
y_{v_{i}}=x_{v_{\eta}}\ ~~~& \ \tau\leq i\leq\eta-1;\\
\\
y_{v_{a(i,j)}}=\max\{x_{v_{a(w,1)}}|\, \tau\leq w\leq \eta\}=\mu\ ~~~& \ \tau\leq i\leq \eta-1, 1\leq j\leq k-2;\\
\\
y_{v_{a(\tau-1,j)}}=\max\{x_{v_{a(\tau-1,1)}}, \mu\}\ ~~~& \ 1\leq j\leq k-2;\\
\\
y_{v_{i}}=x_{v_{i}}\ ~~~& \ 1\leq i\leq\tau-1;\\
\\
y_{v_{i}}=x_{v_{i}}\ ~~~& \ \eta+1\leq i\leq\frac{g+1}{2};\\
\\
y_{v_{g-i+1}}=y_{v_{1+i}}\ ~~~& \ 1\leq i\leq \frac{g-1}{2}\\
\\
y_{v}=x_{v}\ ~~~& \ others. \end{array}\right.$$

Note that $$(\mathcal {A}_{\mathbb{G}_{1,g}}Y)_{v_{\eta}}=\frac{y_{v_{\eta+1}}+y_{v_{\eta-1}}
+(k-2)y_{v_{a(\eta,1)}}+(k-2)y_{v_{a(\eta-1,1)}}}{k-1}\hspace{3cm}$$
$$=\frac{x_{v_{\eta+1}}+x_{v_{\eta}}
+(k-2)x_{v_{a(\eta,1)}}+(k-2)y_{v_{a(\eta-1,1)}}}{k-1}$$
$$>\frac{x_{v_{\eta+1}}+x_{v_{\eta-1}}
+(k-2)x_{v_{a(\eta,1)}}+(k-2)x_{v_{a(\eta-1,1)}}}{k-1}$$
$$=\rho(\mathbb{G}_{1,g})x_{v_{\eta}}=\rho(\mathbb{G}_{1,g})y_{v_{\eta}}.$$

Without loss of generality, we suppose $\mu=x_{v_{a(\eta-1,1)}}$. Then for $1\leq j\leq k-2$,
$$(\mathcal {A}_{\mathbb{G}_{1,g}}Y)_{v_{a(\eta-1,j)}}=\frac{y_{v_{\eta}}+y_{v_{\eta-1}}
+(k-3)y_{v_{a(\eta-1,1)}}}{k-1}\hspace{3cm}$$
$$=\frac{x_{v_{\eta}}+x_{v_{\eta}}
+(k-3)x_{v_{a(\eta-1,1)}}}{k-1}$$
$$>\frac{x_{v_{\eta}}+x_{v_{\eta-1}}
+(k-3)x_{v_{a(\eta-1,1)}}}{k-1}$$
$$=\rho(\mathbb{G}_{1,g})x_{v_{a(\eta-1,1)}}=\rho(\mathbb{G}_{1,g})y_{v_{a(\eta-1,1)}}.$$
In the same way, for $\tau\leq i\leq\eta-2$, $1\leq j\leq k-2$, we get that
$$(\mathcal {A}_{\mathbb{G}_{1,g}}Y)_{v_{a(i,j)}}\geq \rho(\mathbb{G}_{1,g})y_{v_{a(i,j)}};$$
for $\tau\leq i\leq\eta-1$, we get that $$(\mathcal {A}_{\mathbb{G}_{1,g}}Y)_{v_{i}}\geq \rho(\mathbb{G}_{1,g})y_{v_{i}};$$
for other vertices, we get that $$(\mathcal {A}_{\mathbb{G}_{1,g}}Y)_{v}\geq \rho(\mathbb{G}_{1,g})y_{v}.$$
Then using Lemma \ref{le3,04,01} gets a contradiction that $\rho(\mathbb{G}_{1,g})>\rho(\mathbb{G}_{1,g})$.

{\bf Case 2} $x_{v_{\eta}}> x_{v_{i}}$ for $1\leq i\leq \eta-1$.

Let $Y$ be
a vector on $\mathbb{G}_{1,g}$ satisfying that
$$\left \{\begin{array}{ll}
y_{v_{i}}=x_{v_{\eta}}\ ~~~& \ 1\leq i\leq\eta-1;\\
\\
y_{v_{a(i,j)}}=\max\{x_{v_{a(w,1)}}|\, 1\leq w\leq \eta\}=\mu\ ~~~& \ 1\leq i\leq \eta-1, 1\leq j\leq k-2;\\
\\
y_{v_{i}}=x_{v_{i}}\ ~~~& \ \eta+1\leq i\leq\frac{g+1}{2};\\
\\
y_{v_{g-i+1}}=y_{v_{1+i}}\ ~~~& \ 1\leq i\leq \frac{g-1}{2}\\
\\
y_{v}=x_{v}\ ~~~& \ others. \end{array}\right.$$

In the same way as Case 1, we get a a contradiction that $\rho(\mathbb{G}_{1,g})>\rho(\mathbb{G}_{1,g})$.

From the above two cases, it follows that $x_{v_{i}}\geq x_{v_{i+1}}$ for $i=1, 2, \ldots, \frac{g-1}{2}$. Thus (1.1) follows.
(1.2) is proved similarly.

Subsequently, we prove (2).

(2) By Lemma \ref{le03,01}, we know that $x_{v_{a(1,j)}}=x_{v_{a(1,2)}}$ for $j=3, 4, \ldots, k-2$; $x_{v_{a(i,j)}}=x_{v_{a(i,1)}}$ for $i=2, 3, \ldots, g$, $j=1, 2, \ldots, k-2$.

(2.1) As (1.1), we get that $x_{v_{1}}=x_{v_{2}}$, $x_{v_{2+i}}=x_{v_{g-i+1}}$ for $i=1, 2, \ldots, \frac{g-3}{2}$, $x_{v_{a(1+i,1)}}=x_{v_{a(g-i+1,1)}}$ for $i=1, 2, \ldots, \frac{g-1}{2}$. Next, we prove $x_{v_{i}}\geq x_{v_{i+1}}$ for $i=2, 3, \ldots, \frac{g+1}{2}$.  We prove it also by contradiction. Suppose that there exists some $3\leq i \leq\frac{g+3}{2}$ such that $x_{v_{i}}> x_{v_{i-1}}$. Along the hyperpath $P_{2}=v_{\frac{g+3}{2}}e_{\frac{g+1}{2}}v_{\frac{g+1}{2}}e_{\frac{g-1}{2}}v_{\frac{g-1}{2}}\cdots e_{2}v_{2}$, suppose $\eta$ is the largest number form $3$ to $\frac{g+3}{2}$ such that $x_{v_{\eta}}> x_{v_{\eta-1}}$.

{\bf Case 1} There is some number $3\leq\tau<\eta\leq\frac{g+3}{2}$ that $x_{v_{\eta}}> x_{v_{i}}$ for $\tau\leq i\leq\eta-1$, but $x_{v_{\eta}}\leq x_{v_{\tau-1}}$. In the same way as Case 1 in the proof of (1.1), we get a a contradiction that $\rho(\mathbb{G}_{2,g})>\rho(\mathbb{G}_{2,g})$.

{\bf Case 2} $x_{v_{\eta}}> x_{v_{i}}$ for $2\leq i\leq \eta-1$. Now, we need a claim.

{\bf Claim} $x_{v_{a(1,1)}}>x_{v_{a(1,j)}}$ for $2\leq j\leq k-2$. Note that $x_{v_{a(1,j)}}=x_{v_{a(1,2)}}$ for $j=3, 4, \ldots, k-2$.
Now, it is enough to prove that $x_{v_{a(1,1)}}>x_{v_{a(1,2)}}$. Otherwise, suppose that $x_{v_{a(1,1)}}\leq x_{v_{a(1,2)}}$. Then we let $e^{'}_{\varepsilon}=(e_{\varepsilon}\setminus\{v_{a(1,1)}\})\cup\{v_{a(1,2)}\}$ for $e_{\varepsilon}\in E(G_{2})$ with $v_{a(1,1)}\in e_{\varepsilon}$, $$\mathbb{G}^{'}_{2,g}=\mathbb{G}_{2,g}-\sum_{e_{\varepsilon}\in E(G_{2}), v_{a(1,1)}\in e_{\varepsilon}}e_{\varepsilon}+\sum_{e_{\varepsilon}\in E(G_{2}), v_{a(1,1)}\in e_{\varepsilon}}e^{'}_{\varepsilon}.$$ Using Lemma \ref{le03.02} gets that $\rho(\mathbb{G}^{'}_{2,g})>\rho(\mathbb{G}_{2,g})$, which causes a contradiction because $\mathbb{G}^{'}_{2,g}\cong \mathbb{G}_{2,g}$. Thus our claim holds.

Let $Y$ be
a vector on $\mathbb{G}_{2,g}$ satisfying that
$$\left \{\begin{array}{ll}
y_{v_{i}}=x_{v_{\eta}}\ ~~~& \ 2\leq i\leq\eta-1;\\
\\
y_{v_{a(i,j)}}=\max\{x_{v_{a(w,1)}}|\, 2\leq w\leq \eta\}=\mu\ ~~~& \ 2\leq i\leq \eta-1, 1\leq j\leq k-2;\\
\\
y_{v_{a(1,j)}}=\max\{x_{v_{a(1,2)}}, \mu\}\ ~~~& \ 2\leq j\leq k-2;\\
\\
y_{v_{i}}=x_{v_{i}}\ ~~~& \ \eta+1\leq i\leq\frac{g+3}{2};\\
\\
y_{v_{1}}=y_{v_{2}}\\
\\
y_{v_{g-i+1}}=y_{v_{2+i}}\ ~~~& \ i=1, 2, \ldots, \frac{g-3}{2};\\
\\
y_{v_{a(g-i+1,j)}}=y_{v_{a(1+i,j)}}\ ~~~& \ i=1, 2, \ldots, \frac{g-1}{2}, 1\leq j\leq k-2;\\
\\
y_{v}=\frac{y_{v_{1}}+y_{v_{2}}+(k-3)y_{v_{a(1,2)}}}{x_{v_{1}}+x_{v_{2}}+(k-3)x_{v_{a(1,2)}}}x_{v}\ ~~~& \ v\in V(G_{2}). \end{array}\right.$$

In the same way as Case 1 in the proof of (1.1), we get a a contradiction that $\rho(\mathbb{G}_{2,g})>\rho(\mathbb{G}_{2,g})$.

From the above two cases, it follows that $x_{v_{i}}\geq x_{v_{i+1}}$ for $i=2, 3, \ldots, \frac{g+1}{2}$. Thus (2.1) follows. (2.2) is proved similarly.

Thus the result follows. This completes the proof. \ \ \ \ \ $\Box$
\end{proof}

\begin{lemma} \label{le3,15}
In $\mathbb{G}_{1,g}$ and $G_{2}$, if $|E(G_{1})|=|E(G_{2})|\geq 1$, then $\rho(\mathbb{G}_{2,g})<\rho(\mathbb{G}_{1,g})$.
\end{lemma}

\begin{proof}
Let $X$ be
the principal eigenvector of $\mathbb{G}_{2,g}$. If $x_{v_{1}}\geq x_{v_{a(1,1)}}$, we let $e^{'}_{\varepsilon}=(e_{\varepsilon}\setminus\{v_{a(1,1)}\})\cup\{v_{1}\}$ for $e_{\varepsilon}\in E(G_{2})$ with $v_{a(1,1)}\in e_{\varepsilon}$, $$\mathbb{G}^{'}_{2,g}=\mathbb{G}_{2,g}-\sum_{e_{\varepsilon}\in E(G_{2}), v_{a(1,1)}\in e_{\varepsilon}}e_{\varepsilon}+\sum_{e_{\varepsilon}\in E(G_{2}), v_{a(1,1)}\in e_{\varepsilon}}e^{'}_{\varepsilon};$$ if $x_{v_{1}}< x_{v_{a(1,1)}}$, we let $e^{'}_{g}=(e_{g}\setminus\{v_{1}\})\cup \{v_{a(1,1)}\}$, and $\mathbb{G}^{'}_{2,g}=\mathbb{G}_{2,g}-e_{g}+e^{'}_{g}$. Using Lemma \ref{le03.02} gets that $\rho(\mathbb{G}^{'}_{2,g})>\rho(\mathbb{G}_{2,g})$. Note that $\mathbb{G}^{'}_{2,g}\cong \mathbb{G}_{1,g}$. Consequently, it follows that $\rho(\mathbb{G}_{2,g})<\rho(\mathbb{G}_{1,g})$.
Thus the result follows. This completes the proof. \ \ \ \ \ $\Box$
\end{proof}

Note that $Lop_{g,k,s;1}$ is a special $\mathbb{G}_{1,g}$; $Lop_{g,k,s;2}$ is a special $\mathbb{G}_{2,g}$. Then we have the following Corollary \ref{le3,15,00}.

\begin{corollary} \label{le3,15,00}
If $s\geq 1$, then $\rho(Lop_{g,k,s;2})<\rho(Lop_{g,k,s;1})$.
\end{corollary}

\setlength{\unitlength}{0.7pt}
\begin{center}
\begin{picture}(538,132)
\put(0,129){\circle*{4}}
\put(62,80){\circle*{4}}
\qbezier(0,129)(46,129)(62,80)
\qbezier(0,129)(12,83)(62,80)
\put(2,36){\circle*{4}}
\qbezier(62,80)(12,76)(2,36)
\qbezier(62,80)(42,35)(2,36)
\put(142,79){\circle*{4}}
\qbezier(62,80)(103,104)(142,79)
\qbezier(62,80)(105,58)(142,79)
\put(199,128){\circle*{4}}
\qbezier(142,79)(153,125)(199,128)
\qbezier(142,79)(188,90)(199,128)
\put(199,34){\circle*{4}}
\qbezier(142,79)(155,35)(199,34)
\qbezier(142,79)(186,76)(199,34)
\put(62,87){$u$}
\put(136,86){$v$}
\put(2,91){\circle*{4}}
\put(2,80){\circle*{4}}
\put(2,69){\circle*{4}}
\put(198,91){\circle*{4}}
\put(198,80){\circle*{4}}
\put(198,69){\circle*{4}}
\put(306,95){\circle*{4}}
\put(306,84){\circle*{4}}
\put(306,73){\circle*{4}}
\put(99,78){$e$}
\put(299,130){\circle*{4}}
\put(345,85){\circle*{4}}
\qbezier(299,130)(336,128)(345,85)
\qbezier(299,130)(301,95)(345,85)
\put(301,33){\circle*{4}}
\qbezier(345,85)(343,42)(301,33)
\qbezier(345,85)(305,75)(301,33)
\put(412,85){\circle*{4}}
\qbezier(345,85)(380,107)(412,85)
\qbezier(345,85)(383,63)(412,85)
\put(477,85){\circle*{4}}
\qbezier(412,85)(448,110)(477,85)
\qbezier(412,85)(447,64)(477,85)
\put(518,132){\circle*{4}}
\qbezier(477,85)(484,130)(518,132)
\put(514,34){\circle*{4}}
\qbezier(477,85)(477,44)(514,34)
\qbezier(518,132)(515,99)(477,85)
\qbezier(477,85)(516,71)(514,34)
\put(346,92){$u$}
\put(470,92){$v$}
\put(373,83){$e_{1}$}
\put(442,85){$e_{2}$}
\put(408,91){$w$}
\put(517,94){\circle*{4}}
\put(517,83){\circle*{4}}
\put(517,72){\circle*{4}}
\put(95,25){$D_{1}$}
\put(411,26){$D_{2}$}
\put(190,-9){Fig. 3.4.  $D_{1}$, $D_{2}$}
\end{picture}
\end{center}

Let $e=\{u, z_{1}, z_{2}, \ldots, z_{k-2}, v\}$ be an edge in hypergraph $D_{1}$ with $deg_{D_{1}}(z_{i})=1$ for $i=1, 2, \ldots, k-2$. Let $e_{1}=\{u, z_{1}, z_{2}, \ldots, z_{k-2}, w\}$, $e_{2}=\{w, a_{1}, a_{2}, \ldots, a_{k-2}, v\}$, $D_{2}=D_{1}-e+e_{1}+e_{2}$. We call $D_{2}$ is obtained from $D_{1}$ by subdividing edge $e$, or say that $D_{2}$ is a subdivision of $D_{1}$ by subdividing edge $e$. Let $\mathbb{G}_{1,g+1}$ be a subdivision of $\mathbb{G}_{1,g}$ by subdividing edge $e_{\frac{g+1}{2}}$ if $g$ is odd; $\mathbb{G}_{1,g+1}$ be a subdivision of $\mathbb{G}_{1,g}$ by subdividing edge $e_{\frac{g}{2}}$ if $g$ is even. Let $\mathbb{G}_{2,g+1}$ be a subdivision of $\mathbb{G}_{2,g}$ by subdividing edge $e_{\frac{g+1}{2}}$ if $g$ is odd; $\mathbb{G}_{2,g+1}$ be a subdivision of $\mathbb{G}_{2,g}$ by subdividing edge $e_{\frac{g+2}{2}}$ if $g$ is even.

\begin{lemma} \label{le3,16}
(1) If $|E(G_{1})|\geq 1$ in both $\mathbb{G}_{1,g+1}$ and $\mathbb{G}_{1,g}$, then $\rho(\mathbb{G}_{1,g+1})<\rho(\mathbb{G}_{1,g})$.

(2) If $|E(G_{2})|\geq 1$ in both $\mathbb{G}_{2,g+1}$ and $\mathbb{G}_{2,g}$, then $\rho(\mathbb{G}_{2,g+1})<\rho(\mathbb{G}_{2,g})$.
\end{lemma}

\begin{proof}
(1) Let $X$ be
the principal eigenvector of $\mathbb{G}_{1,g}$. By Lemma \ref{le3,19}, we know that if $g$ is odd, $x_{v_{\frac{g+1}{2}}}=\min\{x_{v_{i}}| 1\leq i\leq g\}$; if $g$ is even, $x_{v_{\frac{g+2}{2}}}=\min\{x_{v_{i}}| 1\leq i\leq g\}$. Using Lemma \ref{le3,13} gets that $\rho(\mathbb{G}_{1,g+1})<\rho(\mathbb{G}_{1,g})$.

(2) is proved as (1).
Thus the result follows. This completes the proof. \ \ \ \ \ $\Box$
\end{proof}

\begin{Prof}
This Theorem follows from the fact that $Lop_{g,k,s;1}$ of order $n$ is a special $\mathbb{G}_{1,g}$, $Lop_{g,k,s;2}$ of order $n$ is a special $\mathbb{G}_{2,g}$, and Lemma \ref{le3,16}.
This completes the proof. \ \ \ \ \ $\Box$
\end{Prof}

\begin{lemma} \label{le3,20}
(1) $\rho(Lop_{g+1,k,s;1})<\rho(Lop_{g+1,k,s-1;1})$ where $s\geq 1$.

(2) $\rho(Lop_{g+1,k,s;2})<\rho(Lop_{g+1,k,s-1;2})$ where $s\geq 1$.
\end{lemma}

\begin{proof}
This Lemma follows from Lemma \ref{le3,14} and Theorem \ref{th01.07}. \ \ \ \ \ $\Box$
\end{proof}

\begin{Proof}
(1) Denote by $\mathcal{C}=v_{1}e_{1}v_{2}e_{2}\cdots v_{g-1}e_{g-1}v_{g}e_{g}v_{1}$ the unique hypercycle in $\mathcal{G}$. By Lemma \ref{le3,14}, it follows that $\rho(\mathcal{C})<\rho(\mathcal{G})$. Then combining Corollary \ref{cl01,04} gets the result.

(2) Denote by $\mathcal{C}=v_{1}e_{1}v_{2}e_{2}\cdots v_{g-1}e_{g-1}v_{g}e_{g}v_{1}$ the unique hypercycle in $\mathcal{G}$ for some $g\leq \frac{n-k+1}{k-1}$.
Note that $\mathcal{G}$ is connected and $\mathcal{G}\ncong \mathcal{C}(n, k)$. Then there is a subgraph in $\mathcal{G}$ is isomorphic to $Lop_{g,k,1;1}$ or isomorphic to $Lop_{g,k,1;2}$. Using Lemma \ref{le3,14}, Theorem \ref{th01.07} and Lemma \ref{le3,20} repeatedly gets the result.
This completes the proof. \ \ \ \ \ $\Box$
\end{Proof}

\begin{Proo}
This Lemma follows from using Lemma \ref{le03.02} repeatedly.
This completes the proof. \ \ \ \ \ $\Box$
\end{Proo}

Denote by $\mathcal{G}(\mathcal{D}v; p, q; v_{p+q}\mathcal{H})$ the $k$-uniform connected hypergraph obtained from $k$-uniform hypergraph $\mathcal{D}$ and $k$-uniform hypergraph $\mathcal{H}$ by adding a pendant path $P_{1}$ with length $p$ at vertex $v$ of $\mathcal{D}$, and adding a path $P_{2}$ with length $q$ between vertex $v$ and vertex $v_{p+q}$ of $\mathcal{H}$, where $\mathcal{D}$ and $\mathcal{H}$ are two disjoint, $V(P_{1})\cap V(\mathcal{D})=\{v\}$, $V(P_{2})\cap V(\mathcal{D})=\{v\}$, $V(P_{2})\cap V(\mathcal{H})=\{v_{p+q}\}$ (see two examples in Fig. 3.5). In particular, if $H=v_{p+q}$, we denote by $\mathcal{G}(\mathcal{D}v; p, q; v_{p+q})$ for $\mathcal{G}(\mathcal{D}v; p, q; v_{p+q}H)$ for short.

\

\setlength{\unitlength}{0.6pt}
\begin{center}
\begin{picture}(765,153)
\qbezier(0,97)(0,108)(9,116)\qbezier(9,116)(19,124)(33,124)\qbezier(33,124)(46,124)(56,116)\qbezier(56,116)(66,108)(66,97)
\qbezier(66,97)(66,85)(56,77)\qbezier(56,77)(46,70)(33,70)\qbezier(33,70)(19,70)(9,77)\qbezier(9,77)(0,85)(0,97)
\put(66,98){\circle*{4}}
\put(110,130){\circle*{4}}
\put(133,130){\circle*{4}}
\put(143,130){\circle*{4}}
\put(123,130){\circle*{4}}
\put(155,130){\circle*{4}}
\put(210,130){\circle*{4}}
\qbezier(211,131)(211,140)(219,146)\qbezier(219,146)(227,153)(239,153)\qbezier(239,153)(250,153)(258,146)\qbezier(258,146)(267,140)(267,131)
\qbezier(267,131)(267,121)(258,115)\qbezier(258,115)(250,109)(239,109)\qbezier(239,109)(227,109)(219,115)\qbezier(219,115)(211,121)(211,131)
\put(110,61){\circle*{4}}
\put(184,61){\circle*{4}}
\put(194,61){\circle*{4}}
\put(174,61){\circle*{4}}
\put(204,61){\circle*{4}}
\put(257,61){\circle*{4}}
\put(18,92){$\mathcal{D}$}
\put(55,94){$v$}
\put(107,68){$v_{1}$}
\put(260,58){$v_{p}$}
\put(95,138){$v_{p+1}$}
\put(214,126){$v_{p+q}$}
\put(250,127){$\mathcal{H}$}
\put(71,20){$\mathcal{G}(Dv; p, q; v_{p+q}\mathcal{H})$}
\qbezier(306,96)(306,107)(315,115)\qbezier(315,115)(325,123)(339,123)\qbezier(339,123)(352,123)(362,115)\qbezier(362,115)(372,107)(372,96)
\qbezier(372,96)(372,84)(362,76)\qbezier(362,76)(352,69)(339,69)\qbezier(339,69)(325,69)(315,76)\qbezier(315,76)(306,84)(306,96)
\put(372,95){\circle*{4}}
\put(420,95){\circle*{4}}
\put(479,95){\circle*{4}}
\put(491,95){\circle*{4}}
\put(502,95){\circle*{4}}
\put(513,96){\circle*{4}}
\put(561,96){\circle*{4}}
\put(609,96){\circle*{4}}
\put(621,96){\circle*{4}}
\put(633,96){\circle*{4}}
\put(644,96){\circle*{4}}
\put(657,96){\circle*{4}}
\put(709,96){\circle*{4}}
\qbezier(709,100)(709,109)(717,116)\qbezier(717,116)(725,124)(737,124)\qbezier(737,124)(748,124)(756,116)\qbezier(756,116)(765,109)(765,100)
\qbezier(765,100)(765,90)(756,83)\qbezier(756,83)(748,76)(737,76)\qbezier(737,76)(725,76)(717,83)\qbezier(717,83)(709,90)(709,100)
\put(471,21){$\mathcal{G}(\mathcal{D}v; 0, p+q; v_{p+q}\mathcal{H})$}
\put(205,-9){Fig. 3.5. $\mathcal{G}(Dv; p, q; v_{p+q}\mathcal{H})$ and $\mathcal{G}(\mathcal{D}v; 0, p+q; v_{p+q}\mathcal{H})$}
\put(329,92){$\mathcal{D}$}
\put(746,96){$\mathcal{H}$}
\put(361,92){$v$}
\put(553,105){$v_{p}$}
\put(600,104){$v_{p+1}$}
\put(712,93){$v_{p+q}$}
\put(413,102){$v_{1}$}
\put(387,80){$e_{1}$}
\put(466,95){\circle*{4}}
\put(436,78){$e_{2}$}
\put(459,102){$v_{2}$}
\put(531,80){$e_{p}$}
\put(162,61){\circle*{4}}
\put(156,66){$v_{2}$}
\put(68,66){$e_{1}$}
\put(126,45){$e_{2}$}
\put(224,45){$e_{p}$}
\put(55,125){$e_{p+1}$}
\put(572,82){$e_{p+1}$}
\qbezier(66,98)(81,126)(110,130)
\qbezier(110,130)(99,105)(66,98)
\qbezier(155,130)(181,143)(210,130)
\qbezier(155,130)(180,119)(210,130)
\qbezier(66,98)(91,89)(110,61)
\qbezier(66,98)(76,69)(110,61)
\qbezier(110,61)(135,73)(162,61)
\qbezier(110,61)(135,52)(162,61)
\qbezier(204,61)(231,72)(257,61)
\qbezier(204,61)(231,52)(257,61)
\qbezier(372,95)(396,105)(420,95)
\qbezier(372,95)(396,86)(420,95)
\qbezier(420,95)(444,106)(466,95)
\qbezier(420,95)(444,87)(466,95)
\qbezier(513,96)(538,108)(561,96)
\qbezier(513,96)(539,89)(561,96)
\qbezier(561,96)(585,109)(609,96)
\qbezier(561,96)(586,88)(609,96)
\qbezier(657,96)(683,109)(709,96)
\qbezier(657,96)(685,89)(709,96)
\end{picture}
\end{center}

\begin{lemma}{\bf \cite{GYSLZ}} \label{le3,39}
If $p, q> 0$, then
$\rho(\mathcal{G}(\mathcal{D}v; p, q; v_{p+q}\mathcal{H})) > \rho(\mathcal{G}(\mathcal{D}v; 0, p+q; v_{p+q}\mathcal{H}))$ (see Fig. 3.5).
\end{lemma}

\

\setlength{\unitlength}{0.7pt}
\begin{center}
\begin{picture}(533,287)
\qbezier(112,76)(112,87)(121,95)\qbezier(121,95)(131,103)(145,103)\qbezier(145,103)(158,103)(168,95)\qbezier(168,95)(178,87)(178,76)
\qbezier(178,76)(178,64)(168,56)\qbezier(168,56)(158,49)(145,49)\qbezier(145,49)(131,49)(121,56)\qbezier(121,56)(112,64)(112,76)
\put(177,76){\circle*{4}}
\put(289,75){\circle*{4}}
\put(363,73){\circle*{4}}
\put(375,73){\circle*{4}}
\put(386,73){\circle*{4}}
\put(216,-9){Fig. 3.6. $\mathcal{G}_{1}$, $\mathcal{G}_{2}$}
\put(129,74){$\mathcal{D}$}
\put(161,74){$v_{0}$}
\put(285,63){$u_{1}$}
\put(400,74){\circle*{4}}
\put(463,74){\circle*{4}}
\put(526,75){\circle*{4}}
\put(451,83){$u_{s-1}$}
\put(531,72){$u_{s}$}
\put(231,54){$\tilde{e}_{1}$}
\put(422,54){$\tilde{e}_{s-1}$}
\put(495,54){$\tilde{e}_{s}$}
\qbezier(177,76)(232,98)(289,75)
\qbezier(177,76)(224,58)(289,75)
\put(202,131){\circle*{4}}
\put(202,78){\circle*{4}}
\qbezier(202,131)(187,102)(202,78)
\qbezier(202,131)(216,111)(202,78)
\put(158,155){\circle*{4}}
\qbezier(202,131)(176,133)(158,155)
\qbezier(158,155)(197,155)(202,131)
\put(250,135){\circle*{4}}
\put(231,77){\circle*{4}}
\qbezier(250,135)(227,115)(231,77)
\qbezier(250,135)(254,96)(231,77)
\put(285,181){\circle*{4}}
\qbezier(250,135)(254,167)(285,181)
\qbezier(285,181)(283,150)(250,135)
\put(392,182){\circle*{4}}
\put(340,182){\circle*{4}}
\qbezier(392,182)(366,198)(340,182)
\qbezier(392,182)(368,169)(340,182)
\put(263,77){\circle*{4}}
\put(295,115){\circle*{4}}
\qbezier(263,77)(265,105)(295,115)
\qbezier(295,115)(295,92)(263,77)
\put(344,132){\circle*{4}}
\qbezier(295,115)(315,139)(344,132)
\qbezier(295,115)(324,111)(344,132)
\put(380,132){\circle*{4}}
\put(369,132){\circle*{4}}
\put(357,132){\circle*{4}}
\put(395,132){\circle*{4}}
\put(449,132){\circle*{4}}
\qbezier(395,132)(421,148)(449,132)
\qbezier(395,132)(421,121)(449,132)
\put(323,181){\circle*{4}}
\put(312,181){\circle*{4}}
\put(300,181){\circle*{4}}
\put(145,155){\circle*{4}}
\put(134,155){\circle*{4}}
\put(122,155){\circle*{4}}
\put(57,154){\circle*{4}}
\put(109,155){\circle*{4}}
\qbezier(57,154)(85,168)(109,155)
\qbezier(57,154)(85,139)(109,155)
\qbezier(66,259)(66,268)(75,274)\qbezier(75,274)(84,281)(97,281)\qbezier(97,281)(109,281)(118,274)\qbezier(118,274)(128,268)(128,259)
\qbezier(128,259)(128,249)(118,243)\qbezier(118,243)(109,237)(97,237)\qbezier(97,237)(84,237)(75,243)\qbezier(75,243)(66,249)(66,259)
\put(127,258){\circle*{4}}
\put(186,258){\circle*{4}}
\qbezier(127,258)(158,274)(186,258)
\qbezier(127,258)(157,247)(186,258)
\put(237,258){\circle*{4}}
\qbezier(186,258)(210,273)(237,258)
\qbezier(186,258)(213,248)(237,258)
\put(276,257){\circle*{4}}
\put(265,257){\circle*{4}}
\put(253,257){\circle*{4}}
\put(290,258){\circle*{4}}
\put(345,258){\circle*{4}}
\qbezier(290,258)(316,271)(345,258)
\qbezier(290,258)(319,247)(345,258)
\put(402,258){\circle*{4}}
\qbezier(345,258)(376,271)(402,258)
\qbezier(345,258)(375,246)(402,258)
\put(442,258){\circle*{4}}
\put(431,258){\circle*{4}}
\put(420,258){\circle*{4}}
\put(456,258){\circle*{4}}
\put(506,258){\circle*{4}}
\qbezier(456,258)(481,271)(506,258)
\qbezier(456,258)(482,247)(506,258)
\put(511,256){$v_{t}$}
\put(340,266){$v_{p}$}
\put(82,256){$\mathcal{D}$}
\put(258,220){$\mathcal{G}_{1}$}
\put(263,15){$\mathcal{G}_{2}$}
\put(121,137){$\mathbf{P}_{1}$}
\put(312,165){$\mathbf{P}_{2}$}
\put(363,111){$\mathbf{P}_{f}$}
\put(351,74){\circle*{4}}
\qbezier(289,75)(323,89)(351,74)
\qbezier(289,75)(323,61)(351,74)
\put(315,55){$\tilde{e}_{2}$}
\put(362,39){$\mathbf{P}_{0}$}
\put(112,256){$v_{0}$}
\put(150,242){$e_{1}$}
\put(181,265){$v_{1}$}
\put(234,264){$v_{2}$}
\put(204,242){$e_{2}$}
\put(477,242){$e_{t}$}
\put(274,278){$\mathbf{P}$}
\put(392,266){$v_{p+1}$}
\qbezier(400,74)(430,89)(463,74)
\qbezier(400,74)(432,61)(463,74)
\qbezier(463,74)(500,90)(526,75)
\qbezier(463,74)(502,62)(526,75)
\put(361,162){\circle*{4}}
\put(368,155){\circle*{4}}
\put(375,148){\circle*{4}}
\end{picture}
\end{center}

\begin{lemma}{\bf \cite{GYSLZ}} \label{le3,41}
$\mathcal{D}$ is a $k$-uniform hypergraph where $v_{0}\in V(\mathcal{D})$. Both $\mathbf{P}=v_{0}e_{1}v_{1}e_{2}v_{2}\cdots e_{t}v_{t}$ and
$\mathbf{P}_{0}=v_{0}\tilde{e}_{1}u_{1}\tilde{e}_{2}u_{2}\cdots \tilde{e}_{s}u_{s}$ are $k$-uniform hyperpaths where $\tilde{e}_{1}=\{v_{0}$, $v_{\varphi(1,1)}$, $v_{\varphi(1,2)}$, $\ldots$, $v_{\varphi(1,k-2)}$, $u_{1}\}$, $V(\mathcal{D})\cap V(\mathbf{P})=\{v_{0}\}$, $V(\mathcal{D})\cap V(\mathbf{P}_{0})=\{v_{0}\}$.
$\mathbf{P}_{1}$, $\mathbf{P}_{2}$, $\ldots$, $\mathbf{P}_{f}$ ($1\leq f\leq k-2$) are $k$-uniform hyperpaths attached respectively at vertices $v_{\varphi(1,1)}$, $v_{\varphi(1,2)}$, $\ldots$, $v_{\varphi(1,f)}$ in $\tilde{e}_{1}$ satisfying $1\leq L(\mathbf{P}_{1})\leq L(\mathbf{P}_{2})\leq\cdots\leq L(\mathbf{P}_{f})\leq  L(\mathbf{P}_{0})-1$, $\sum^{f}_{i=0}L(\mathbf{P}_{i})=t$, $V(\mathcal{D})\cap V(\mathbf{P}_{i})=\emptyset$ for $1\leq i\leq f$. $K$-uniform hypergraph $\mathcal{G}_{1}$ is a $\mathcal{G}(\mathcal{D}v_{0}; 0, t; v_{t})$ consisting of $\mathcal{D}$ and $\mathbf{P}$;
$k$-uniform hypergraph $\mathcal{G}_{2}$ consists of $\mathcal{D}$ and $\mathbf{P}_{0}$, $\mathbf{P}_{1}$, $\mathbf{P}_{2}$, $\ldots$, $\mathbf{P}_{f}$ (see Fig. 3.6).
Then $\rho(\mathcal{G}_{1})<\rho(\mathcal{G}_{2})$.
\end{lemma}

\begin{Pro}
By using Corollary \ref{le3,15,00}, and using Lemma \ref{le3,39} and Lemma \ref{le3,41} repeatedly, it follows the theorem.
This completes the proof. \ \ \ \ \ $\Box$
\end{Pro}

\small {

}

\end{document}